\documentclass[final,5p,times]{elsarticle}

\usepackage{graphicx,subfigure}
\graphicspath{{./figures/}} 
\usepackage{amsfonts,amsbsy,amscd,amsgen,amsthm,dsfont,amsmath,amssymb,mathrsfs,mathtools,array}
\usepackage[]{natbib}
\usepackage[colorlinks=true,allcolors=blue]{hyperref}
\usepackage{bm,multirow,esvect,cancel}
\theoremstyle{definition}
\usepackage{tabularx}
\usepackage{xcolor,comment,soul}
\usepackage{enumitem}
\usepackage[utf8]{inputenc}
\usepackage[T1]{fontenc}
\usepackage[]{algorithm}
\usepackage{algpseudocode} 
\usepackage{caption}
\DeclareCaptionFormat{algor}{%
	\hrulefill\par\offinterlineskip\vskip1pt%
	\textbf{#1#2}#3\offinterlineskip\hrulefill
}
\DeclareCaptionStyle{algori}{singlelinecheck=off,format=algor,labelsep=space}
\makeatletter
\def\ps@pprintTitle{%
    \let\@oddhead\@empty
    \let\@evenhead\@empty
    \def\@oddfoot{\footnotesize
         {Published in Examples and Counterexamples, \url{https://doi.org/10.1016/j.exco.2022.100077}} \hfill May 3, 2022}%
    \let\@evenfoot\@oddfoot
    }
\makeatother

\newcommand{\id}{\mathrm d}

\DeclareMathAlphabet\mathbfcal{OMS}{cmsy}{b}{n}

\newtheorem{thm}{Theorem}
\newtheorem{rem}{Remark}
\newtheorem{lem}{Lemma}
\newtheorem{defn}{Definition}

\newtheorem{cor}{Corollary}

\journal{Examples and Counterexamples}

\begin{document}

\begin{frontmatter}

\title{Extension of Stein's lemma derived by using an integration by differentiation technique}

\author[1]{Konstantinos Mamis}
\ead{kmamis@ncsu.edu}
\address[1]{Department of Mathematics, North Carolina State University,  2311 Stinson Drive, Raleigh, NC 27695-8205, USA}

\begin{abstract}
We extend Stein's lemma for averages that explicitly contain the Gaussian random variable at a power. We present two proofs for this extension of Stein's lemma, with the first being a rigorous proof by mathematical induction. The alternative, second proof is a constructive formal derivation in which we express the average not as an integral,  but as the action of a pseudodifferential operator defined via the Gaussian moment-generating function.  In extended Stein's lemma, the absolute values of the coefficients of the probabilist's Hermite polynomials appear,  revealing yet another link between Hermite polynomials and normal distribution. 
\end{abstract}



\begin{keyword}
normal distribution \sep Stein's lemma  \sep Hermite polynomials \sep generalized factorial coefficients \sep pseudodifferential operator 
\MSC primary 60E05 \sep secondary 60E10 \sep 47G30 \sep 05A10
\end{keyword}

\end{frontmatter}
\section{Introduction and main results}\label{sec:introduction}
Stein's lemma \cite{Stein1981} is a celebrated result in probability theory with many applications in statistics, see e.g. \cite{Mukhopadhyay2021}.  For a scalar, zero-mean Gaussian random variable $X$ with variance $\sigma^2$, Stein's lemma reads
\begin{equation}\label{eq:Stein}
\mathbb{E}\left[g(X)X\right]=\sigma^2\mathbb{E}\left[g'(X)\right],
\end{equation}
where $\mathbb{E}[\cdot]$ is the mean value operator, and prime denotes the first derivative of function $g$.  By Theorem \ref{thm} of the present work, we extend Stein's lemma for the average $\mathbb{E}\left[g(X)X^n\right]$, which is expressed as a finite series that contains averages of derivatives of $g(X)$ up to the $n$th order. We present two ways of proving Theorem \ref{thm}; in Sec.~\ref{sec:induction}, the extended Stein's lemma is proven rigorously by mathematical induction on index $n$ of power $X^n$ inside the average.  Furthermore, and in order to provide the reader with more insight, we present an alternative, constructive proof in Sec.~\ref{sec:constructive}.  The constructive proof is based on a formal expression of the average as the action of a pseudodifferential operator, introduced by Definition \ref{thm:average} via the Gaussian moment-generating function. Series expansion of this pseudodifferential operator allows us to treat the averages not as integrals, but as Taylor-like infinite series. Thus, and under the formal assumption that all infinite series involved are summable,  we are able to calculate $\mathbb{E}\left[g(X)X^n\right]$ without performing any integrations; instead, we perform only the differentiations that appear in the series terms, justifying thus the name integration by differentiation for this technique.   
\begin{thm}[Extension of Stein's lemma]\label{thm}
Let $X$ be a scalar Gaussian variable with zero mean value and variance $\sigma^2$.  Given an $n\in\mathbb{N}$, let also $g$  be a $C^n(\mathbb{R}\rightarrow\mathbb{R})$ function for which averages $\mathbb{E}\left[g(X)X^n\right]$ and $\mathbb{E}\left[g^{(n-2k)}(X)\right]$ for $k=0,\ldots, \lfloor n/2\rfloor$ exist, with $g^{(\ell)}$ denoting the $\ell$th derivative of $g$,  $g^{(0)}:=g$, and $\lfloor \cdot\rfloor$ being the floor function.  It holds true that:
\begin{equation}\label{eq:ext_Stein}
\mathbb{E}\left[g(X)X^n\right]=\sum_{k=0}^{\lfloor n/2\rfloor}H_{n,k}\sigma^{2(n-k)}\mathbb{E}\left[g^{(n-2k)}(X)\right],
\end{equation}
where
\begin{equation}\label{eq:hermite_numbers}
H_{n,k}=\frac{n!}{2^kk!(n-2k)!}, \ \ k=0,\ldots,\lfloor n/2\rfloor,
\end{equation}
are the absolute values of the coefficients appearing in the $n$th-order \textit{probabilist's Hermite polynomial} $He_{n}(x)$ \cite[expression 22.3.11]{Abramowitz1964}:
\begin{equation}\label{eq:Hermite_polynomials}
He_{n}(x)=\sum_{k=0}^{\lfloor n/2\rfloor}(-1)^k H_{n,k}x^{n-2k}. 
\end{equation}
In the On-Line Encyclopedia of Integer Sequences (OEIS) \cite{OEIS}, $H_{n,k}$ are referred to as the Bessel numbers. In this work, we shall call $H_{n,k}$ the signless Hermite coefficients, as a more suggestive term in our context.
\end{thm}
\begin{proof}
Theorem \ref{thm} is proven by mathematical induction,  performed in Sec.~\ref{sec:induction}.  In Sec.~\ref{sec:constructive}, we also present  a constructive derivation of Eq.~\eqref{eq:ext_Stein} that employs an integration by differentiation technique.
\end{proof}
\begin{rem}[Orders of derivatives in Eq.~\eqref{eq:ext_Stein}]\label{rem:order}
We easily observe that for an even (odd) power $n$ of $X$, only the averages of the even (odd)-order derivatives of $g(X)$, up to the $n$th order,  appear in the right-hand side of Eq.~\eqref{eq:ext_Stein}.  The average of $g(X)$ itself appears in Eq.~\eqref{eq:ext_Stein} only for even powers.  
\end{rem}
Furthermore,  for the complete treatment of $\mathbb{E}\left[g(X)X^n\right]$,  we also have to consider the case where the mean value of $X$ is non-zero. Thus, we state the following Theorem:
\begin{thm}[Extension of Stein's lemma for non-zero mean]\label{thm:mean}
For a Gaussian variable $X$ with non-zero mean value $\mu$, the extended Stein's lemma reads
\begin{equation}\label{eq:non_zero_mu}
\mathbb{E}\left[g(X)X^n\right]=\sum_{\ell=0}^n\binom{n}{\ell}\mu^{n-\ell}\sum_{k=0}^{\lfloor \ell/2\rfloor}H_{\ell,k}\sigma^{2(\ell-k)}\mathbb{E}\left[g^{(\ell-2k)}(X)\right],
\end{equation}
where $\binom{n}{\ell}=\frac{n!}{\ell!(n-\ell)!}$ is the binomial coefficient.
\end{thm}
\begin{proof}
Eq.~\eqref{eq:non_zero_mu} can be proven by mathematical induction, similar to the one of Sec.~\ref{sec:induction} for Theorem \ref{thm}.  Alternatively, one can perform a constructive proof for Eq.~\eqref{eq:non_zero_mu} similar to the one of Sec.~\ref{sec:constructive} for Theorem \ref{thm},  by expressing the averages via Eq.~\eqref{eq:pseudo2}.
\end{proof}
The formulas of extended Stein's lemma, stated in Theorems \ref{thm}, \ref{thm:mean}, are useful in engineering applications where averages $\mathbb{E}\left[g(X)X^n\right]$ arise,  see e.g. the recent work \cite{Skaltsas2021} in tribology. Note that the respective extensions of Stein's lemma for multivariate Gaussian variables are also feasible and will be the topic of a forthcoming work. 
\begin{rem}[The need for extended Stein's lemma]
One could argue that the evaluation of averages $\mathbb{E}\left[g(X)X^n\right]$ can always be performed by repetitive applications of the classical Stein's lemma. However, the advantage of  closed-form formulas \eqref{eq:ext_Stein}, \eqref{eq:non_zero_mu} over the recursive ones is manifested for increasing $n$, see e.g.  \cite[Sec. 2.3]{Skaltsas2021},  where keeping track of the repetitive applications of Eq.~\eqref{eq:Stein} becomes progressively more difficult. On the other hand, the averages and coefficients appearing in the right-hand sides of Eqs.~ \eqref{eq:ext_Stein}, \eqref{eq:non_zero_mu} are a priori known for any $n$, offering thus an easy and tractable calculation. 
\end{rem}
\begin{cor}[Central moments of a Gaussian variable]
The higher order moments of a zero-mean Gaussian variable $X$ are given by the formula
\begin{equation}\label{moment}
\mathbb{E}[X^n]= \left\{
\begin{array}{ll}
      0 & \text{for } n\text{ odd}, \\
      \sigma^n(n-1)!! & \text{for } n\text{ even}.\\
\end{array} 
\right. 
\end{equation} 
\end{cor}
\begin{proof}
Eq.~\eqref{moment} is usually derived by repetitive integrations by parts,  see e.g.  \cite[p. 148]{Papoulis2002}.  Here,  we easily derive it from Eq.~\eqref{eq:ext_Stein}, by considering $g(x)=1$. In this case, all derivatives of $g$ appearing in the right-hand side of Eq.~\eqref{eq:ext_Stein} are zero, except for the zeroth-order one, which equals to 1.  Since  Eq.~\eqref{eq:ext_Stein} does not contain $\mathbb{E}[g(X)]$ for odd $n$ (see Remark \ref{rem:order}),  all odd moments of $X$ are equal to zero.  For even $n=2\ell$, zeroth-order derivative appears for $k=\ell$, and Eq.~\eqref{eq:ext_Stein} for $g(x)=1$ reads
\[\mathbb{E}[X^{2\ell}]=H_{2\ell,\ell}\sigma^{2\ell}=\frac{(2\ell)!}{2^\ell\ell!}\sigma^{2\ell}.\]
Then, by using factorial relations $n!=n!!(n-1)!!$,  and $n!!=2^\ell\ell!$ for $n=2\ell$,  we obtain the branch of Eq.~\eqref{moment} for even $n$.
\end{proof}

\section{Proof of Theorem \ref{thm} by mathematical induction}\label{sec:induction}
Theorem \ref{thm} can be proven  by mathematical induction. For $n=1$, it is easy to see that Eq.~\eqref{eq:ext_Stein} results in Stein's lemma \eqref{eq:Stein}, since, from Eq.~\eqref{eq:hermite_numbers}, $H_{1,0}=1$. We then assume that Eq.~\eqref{eq:ext_Stein} holds true for a particular $n$ (inductive hypothesis).  Now, for $n+1$, and by using the inductive hypothesis, we obtain:
\begin{align}\label{eq1}
\mathbb{E}\left[g(X)X^{n+1}\right]&=\mathbb{E}\left[\big(g(X)X\big)X^n\right]=\nonumber\\&=\sum_{k=0}^{\lfloor n/2\rfloor}H_{n,k}\sigma^{2(n-k)}\mathbb{E}\left[\big(g(X)X\big)^{(n-2k)}\right].
\end{align}
For the evaluation of the derivative inside the average in the right-hand side of Eq.~\eqref{eq1}, we use the general Leibniz rule  \cite[expression 3.3.8]{Abramowitz1964}:
\begin{equation}\label{leibniz}
\big(g(x)x\big)^{(n-2k)}=\sum_{\ell=0}^{n-2k}\binom{n-2k}{\ell}g^{(n-2k-\ell)}(x)x^{(\ell)}.
\end{equation}
Since $x^{(0)}=x$, $x^{(1)}=1$, and $x^{(\ell)}=0$ for $\ell\geq 2$, Eq.~\eqref{leibniz} is simplified into
\begin{equation}\label{leibniz2}
\big(g(x)x\big)^{(n-2k)}=g^{(n-2k)}(x)x+(n-2k)g^{(n-2k-1)}(x),
\end{equation}
under the convention that $g^{(-1)}(x)=0$. By using Leibniz rule~\eqref{leibniz2} in Eq.~\eqref{eq1}, we have
\begin{equation}\label{eq2}
\mathbb{E}\left[g(X)X^{n+1}\right]=\Sigma_1+\Sigma_2,
\end{equation}
where
\begin{equation}\label{S1}
\Sigma_1=\sum_{k=0}^{\lfloor n/2\rfloor}H_{n,k}\sigma^{2(n-k)}\mathbb{E}\left[g^{(n-2k)}(X)X\right],
\end{equation}
and
\begin{equation}\label{S2}
\Sigma_2=\sum_{k=0}^{\lfloor (n-1)/2\rfloor}(n-2k)H_{n,k}\sigma^{2(n-k)}\mathbb{E}\left[g^{(n-2k-1)}(X)\right].
\end{equation}
Note that the upper limit of $k$-sum in the right-hand side of Eq.~\eqref{S2} has been changed from $\lfloor n/2\rfloor$ to $\lfloor (n-1)/2\rfloor$, in order to exclude any term containing $g^{(-1)}$. The only such term is for even $n$ and $k=\lfloor n/2\rfloor$.  Thus,  the upper limit of $k$-sum is $\lfloor n/2\rfloor$ for $n$ odd and $\lfloor n/2\rfloor-1$ for $n$ even.  These two values can be expressed in a unified way as $\lfloor (n-1)/2\rfloor$.  By also performing a change in index,  and using the fact that $\lfloor (n+1)/2\rfloor=\lfloor (n-1)/2\rfloor+1$, $\Sigma_2$ reads
\begin{equation}\label{S2_new}
\Sigma_2=\sum_{k=1}^{\lfloor (n+1)/2\rfloor}(n-2k+2)H_{n,k-1}\sigma^{2(n+1-k)}\mathbb{E}\left[g^{(n+1-2k)}(X)\right].
\end{equation}
For $\Sigma_1$, we apply Stein's lemma~\eqref{eq:Stein} again,  at the average appearing in the right-hand side of Eq.~\eqref{S1}, resulting into
\begin{equation}\label{S1_new}
\Sigma_1=\sum_{k=0}^{\lfloor n/2\rfloor}H_{n,k}\sigma^{2(n+1-k)}\mathbb{E}\left[g^{(n+1-2k)}(X)\right].
\end{equation}
Under Eqs.~\eqref{S2_new}, \eqref{S1_new},  the right-hand side of Eq.~\eqref{eq2} is rearranged into
\begin{align}\label{eq3}
\mathbb{E}\left[g(X)X^{n+1}\right]=\sum_{k=0}^{\lfloor (n+1)/2\rfloor}&\left[(n-2k+2)H_{n,k-1}+H_{n,k}\right]\times\nonumber\\&\times\sigma^{2(n+1-k)}\mathbb{E}\left[g^{(n+1-2k)}(X)\right],
\end{align}
under the convention that $H_{n,k}=0$ for $k<0$ or $k>n/2$. 
\begin{lem}[Recurrence relation for $H_{n,k}$]\label{recur}
For $n\in\mathbb{N}$, $k=0,\ldots,\lfloor(n+1)/2\rfloor$,  and with $H_{n,k}=0$ for $k<0$ or $k>n/2$, it holds true that 
\begin{equation}\label{H}
H_{n+1,k}=(n-2k+2)H_{n,k-1}+H_{n,k}.
\end{equation}
\end{lem}
\begin{proof}
See \ref{A}.  Eq.~\eqref{H} is also stated in the relevant OEIS entry \cite{OEIS}.
\end{proof}
Substitution of Eq.~\eqref{H} into Eq. ~\eqref{eq3} results in
\begin{equation}\label{eq4}
\mathbb{E}\left[g(X)X^{n+1}\right]=\sum_{k=0}^{\lfloor (n+1)/2\rfloor}H_{n+1,k}\sigma^{2(n+1-k)}\mathbb{E}\left[g^{(n+1-2k)}(X)\right],
\end{equation}
which is Eq.~\eqref{eq:ext_Stein} for $n+1$, completing thus the proof of Theorem \ref{thm} by induction.
\section{Formal derivation of Theorem \ref{thm} using integration by differentiation}
\label{sec:constructive}
Stein's lemma \eqref{eq:Stein} is usually proven by integration by parts (see e.g. \cite[lemma 1]{Stein1981}), by employing the definition of mean value operator $\mathbb{E}[\cdot]$ as an integral over $\mathbb{R}$.  That is,  since random variable $X$ follows the univariate normal distribution $f(x)=\left(\sigma\sqrt{2\pi}\right)^{-1}\exp\left(-x^2/(2\sigma^2)\right)$, we calculate:
\begin{align}\label{eq:byparts}
\mathbb{E}\left[g(X)X\right]&=\int_{\mathbb{R}}g(x)xf(x)\id x=-\sigma^2\int_{\mathbb{R}}g(x)f'(x)\id x=\nonumber\\&=\sigma^2\int_{\mathbb{R}}g'(x)f(x)\id x=\sigma^2\mathbb{E}\left[g'(X)\right].
\end{align} 
While one could also go on to try to find the adequate repetitive integrations by parts for the calculation of $\mathbb{E}\left[g(X)X^n\right]$, here we shall follow a different, more tractable path, based on the interpretation of average not as an integral, but as the action of a pseudodifferential operator,  defined as follows.
\begin{defn}[Mean value as the action of an averaged shift operator]\label{thm:average}
Let $X$ be a zero-mean Gaussian random variable with variance $\sigma^2$ and $g$ be a $C^{\infty}(\mathbb{R}\rightarrow\mathbb{R})$ function.  Then, average $\mathbb{E}\left[g(X)\right]$ is expressed as
\begin{equation}\label{eq:pseudo}
\mathbb{E}\left[g(X)\right]=\left.\exp\left(\frac{\sigma^2}{2}\frac{\id^2}{\id x^2}\right)g(x)\right|_{x=0}.
\end{equation}
In Eq.~\eqref{eq:pseudo}, $\exp\left(\left(\sigma^2/2\right)\id^2/\id x^2\right)$ is a pseudodifferential operator, called the averaged shift operator (see \ref{A:average}), whose action is to be understood by its series form 
\begin{equation}\label{eq:series0}
\mathbb{E}\left[g(X)\right]=\sum_{m=0}^\infty\frac{\sigma^{2m}}{2^mm!}\left.\frac{\id^{2m} g(x)}{\id x^{2m}}\right|_{x=0},
\end{equation}
with $\cdot|_{x=0}$ denoting that all derivatives appearing in the right-hand side of Eq.~\eqref{eq:series0} are calculated at $x=0$. 
\end{defn}
\begin{proof}
Definition \ref{thm:average} is formally derived in \ref{A:average} by employing the moment-generating function of $X$. Its infinite dimensional version is presented in \cite{Athanassoulis2019a}, and this concept is also found in \cite[Ch. 4]{Klyatskin2005}.
\end{proof}
\begin{rem}[Integration by differentiation]
Choosing Eq.~\eqref{eq:series0} for the evaluation of $\mathbb{E}\left[g(X)X^n\right]$ is more convenient than integration by parts, since, by using Eq.~\eqref{eq:series0}, only the calculation of derivatives is needed. Other integration by differentiation techniques have also been recently proposed, see e.g. \cite{Jia2017},  exploiting the fact that differentiation is generally easier than integration. 
\end{rem}
Thus, under the additional assumption that $g$ is $C^{\infty}(\mathbb{R}\rightarrow\mathbb{R})$, we express $\mathbb{E}\left[g(X)X^n\right]$ via Eq.~\eqref{eq:series0} as 
\begin{equation}\label{eq:series}
\mathbb{E}\left[g(X)X^n\right]=\sum_{m=0}^\infty\frac{\sigma^{2m}}{2^mm!}\left.\frac{\id^{2m} \left[g(x)x^n\right]}{\id x^{2m}}\right|_{x=0}.
\end{equation}
We evaluate the derivatives appearing in the right-hand side of Eq.~\eqref{eq:series} by using the general Leibniz rule
\begin{equation}\label{eq:Lei}
\frac{\id^{2m} \left[g(x)x^n\right]}{\id x^{2m}}=\sum_{\ell=0}^{2m}\binom{2m}{\ell}g^{(2m-\ell)}(x)\left(x^n\right)^{(\ell)},
\end{equation} 
and since $\left(x^{n}\right)^{(\ell)}=(n!/(n-\ell)!)x^{n-\ell}$ for $n\geq\ell$ and zero for $n<\ell$:
\begin{equation}\label{eq:Lei2}
\frac{\id^{2m} \left[g(x)x^n\right]}{\id x^{2m}}=\sum_{\ell=0}^{\min\{2m,n\}}\binom{2m}{\ell}\frac{n!}{(n-\ell)!}g^{(2m-\ell)}(x)x^{n-\ell}.
\end{equation}
Note that,  in Eq.~\eqref{eq:series},  derivative \eqref{eq:Lei2} has to be calculated for $x=0$.  Since, for $x=0$, all terms of the sum in the right-hand side of \eqref{eq:Lei2} are zero except for $\ell=n$, we obtain:
\begin{equation}\label{eq:der}
\left.\frac{\id^{2m} \left[g(x)x^n\right]}{\id x^{2m}}\right|_{x=0}= \left\{
\begin{array}{ll}
      0 & \text{for }m<n/2, \\
      \frac{(2m)!}{(2m-n)!}\left.g^{(2m-n)}(x)\right|_{x=0} & \text{for }m\geq n/2.\\
\end{array} 
\right. 
\end{equation} 
By substitution of expression ~\eqref{eq:der} into Eq.~\eqref{eq:series}, and after some algebraic manipulations,  we have
\begin{equation}\label{eq:series2}
\mathbb{E}\left[g(X)X^n\right]=\sum_{m=\lceil n/2\rceil}^\infty\frac{\sigma^{2m}}{2^m}\frac{(2m)^{\underline{n}}}{m!}\left.g^{(2m-n)}(x)\right|_{x=0},
\end{equation}
where $\lceil \cdot\rceil$ denotes the ceiling function and $(2m)^{\underline{n}}=(2m)(2m-1)\cdots(2m-n+1)$ is the falling factorial.  In order to evaluate further the right-hand side of Eq.~\eqref{eq:series2}, the falling factorial of $2m$ has to be expressed in terms of the falling factorial of $m$.  Following Charalambides \cite[Sec.  8.4]{Charalambides2002}, this is performed by using the \textit{generalized factorial coefficients with parameter 2} $C(n,\ell;2)$, that have the property
\begin{equation}\label{eq:factorial}
(2m)^{\underline{n}}=\sum_{\ell=0}^nC(n,\ell;2)m^{\underline{\ell}}.
\end{equation}
Generalized factorial coefficients $C(n,\ell;2)$ are defined in terms of the Stirling numbers of first $s(n,k)$ and second $S(n,k)$ kind \cite[theorem 8.13]{Charalambides2002}:
\begin{equation}\label{eq:stirling}
C(n,\ell;2)=\sum_{k=\ell}^n2^\ell s(n,k)S(k,\ell),
\end{equation}
and obey the following recurrence relation for $n\in\mathbb{N}$, $\ell=1,2,\ldots,n+1$ \cite[theorem 8.19]{Charalambides2002}:
\begin{equation}\label{eq:recurenceC}
C(n+1,\ell;2)=(2\ell-n)C(n,\ell;2)+2C(n,\ell-1;2),
\end{equation}
with initial conditions $C(0,0;2)=1$, $C(n,0;2)=0$ for $n>0$, and $C(n,\ell;2)=0$ for $\ell>n$.  Also, since $s(n,n)=S(n,n)=1$,  see \cite[proposition 5.3.2]{Cameron1994}, we deduce from Eq.~\eqref{eq:stirling} that 
\begin{equation}\label{eq:final_cond}
C(n,n;2)=2^n, \ \ \text{for }n\in\mathbb{N}.
\end{equation}
Furthermore,  in \cite[remark 8.8]{Charalambides2002},  it is proved that $C(n,\ell;2)=0$ for $\ell<n/2$.  This property results in Eq.~\eqref{eq:factorial} to be updated into
\begin{equation}\label{eq:factorial2}
(2m)^{\underline{n}}=\sum_{\ell=\lceil n/2\rceil}^nC(n,\ell;2)m^{\underline{\ell}}.
\end{equation}
Last, from the definition of falling factorial, $m^{\underline{n}}$ is zero for $n>m$,  and thus Eq.~\eqref{eq:factorial2} is finally expressed as
\begin{equation}\label{eq:factorial3}
(2m)^{\underline{n}}=\sum_{\ell=\lceil n/2\rceil}^{\min\{m,n\}}C(n,\ell;2)m^{\underline{\ell}}.
\end{equation}
Substituting Eq.~\eqref{eq:factorial3} into Eq.~\eqref{eq:series2}, and use of $m^{\underline{\ell}}/m!=1/(m-\ell)!$ results in
\begin{align}\label{eq:series3}
\mathbb{E}\left[g(X)X^n\right]&=\sum_{m=\lceil n/2\rceil}^\infty\sum_{\ell=\lceil n/2\rceil}^{\min\{m,n\}}\frac{C(n,\ell;2)}{(m-\ell)!}\frac{\sigma^{2m}}{2^m}\left.g^{(2m-n)}(x)\right|_{x=0}=\nonumber\\
&=\sum_{m=\lceil n/2\rceil}^{n}\sum_{\ell=\lceil n/2\rceil}^{m}\frac{C(n,\ell;2)}{(m-\ell)!}\frac{\sigma^{2m}}{2^m}\left.g^{(2m-n)}(x)\right|_{x=0}+\nonumber\\&+\sum_{m=n+1}^\infty\sum_{\ell=\lceil n/2\rceil}^{n}\frac{C(n,\ell;2)}{(m-\ell)!}\frac{\sigma^{2m}}{2^m}\left.g^{(2m-n)}(x)\right|_{x=0}.
\end{align}
Double summations in the rightmost side of Eq.~\eqref{eq:series3} are easily rearranged into
\begin{align}\label{eq:series4}
\mathbb{E}\left[g(X)X^n\right]&=\sum_{\ell=\lceil n/2\rceil}^{n}\sum_{m=\ell}^{n}\frac{C(n,\ell;2)}{(m-\ell)!}\frac{\sigma^{2m}}{2^m}\left.g^{(2m-n)}(x)\right|_{x=0}+\nonumber\\&+\sum_{\ell=\lceil n/2\rceil}^{n}\sum_{m=n+1}^\infty\frac{C(n,\ell;2)}{(m-\ell)!}\frac{\sigma^{2m}}{2^m}\left.g^{(2m-n)}(x)\right|_{x=0}=\nonumber\\
&=\sum_{\ell=\lceil n/2\rceil}^{n}\sum_{m=\ell}^{\infty}\frac{C(n,\ell;2)}{(m-\ell)!}\frac{\sigma^{2m}}{2^m}\left.g^{(2m-n)}(x)\right|_{x=0}.
\end{align}
An index change in the second sum, as well as the use of formula \eqref{eq:series0},  results in
\begin{align}\label{eq:series5}
\mathbb{E}\left[g(X)X^n\right]&=\sum_{\ell=\lceil n/2\rceil}^{n}\frac{C(n,\ell;2)}{2^\ell}\sigma^{2\ell}\sum_{m=0}^\infty\frac{\sigma^{2m}}{2^mm!}\left.g^{(2m+2\ell-n)}(x)\right|_{x=0}=\nonumber\\&=\sum_{\ell=\lceil n/2\rceil}^{n}\frac{C(n,\ell;2)}{2^\ell}\sigma^{2\ell}\mathbb{E}\left[g^{(2\ell-n)}(X)\right].
\end{align}
By also performing the change of summation index $k=n-\ell$ in Eq.~\eqref{eq:series5}, we obtain
\begin{equation}\label{eq:final}
\mathbb{E}\left[g(X)X^n\right]=\sum_{k=0}^{\lfloor n/2\rfloor}\frac{C(n,n-k;2)}{2^{n-k}}\sigma^{2(n-k)}\mathbb{E}\left[g^{(n-2k)}(X)\right].
\end{equation}
Eq.~\eqref{eq:final} is of the same form as the extended Stein's lemma,  Eq.~\eqref{eq:ext_Stein}.  What remains to be proven is the identification of $C(n,n-k;2)/2^{n-k}$ as the signless Hermite coefficients $H_{n,k}$.  This is performed in the following Lemma, which concludes the constructive formal derivation of Theorem \ref{thm}.
\begin{lem}[Signless Hermite coefficients as rearranged, rescaled, generalized factorial coefficients with parameter 2]\label{HC} For $H_{n,k}$ being the signless Hermite coefficients defined by Eq.~\eqref{eq:hermite_numbers} and $C(n,\ell;2)$ being the generalized factorial coefficients with the property \eqref{eq:factorial}, it holds true that
\begin{equation}\label{eq:HC}
H_{n,k}=\frac{C(n,n-k;2)}{2^{n-k}}, \ \ \text{for} \ \ n\in\mathbb{N}, \ \ k=0,\ldots,\lfloor n/2\rfloor.
\end{equation}
\end{lem}
\begin{proof}
See \ref{B}. To the best of our knowledge, relation \eqref{eq:HC} has not been pointed out before. 
\end{proof}

\section*{Declaration of Competing Interest}
The author declares that he has no known competing financial interests or personal relationships that could have appeared to influence the work reported in this paper.

\section*{Acknowledgement}
Part of this work was conducted while the author was a PhD candidate in the National Technical University of Athens,  under the supervision of professor G.A. Athanassoulis.

\appendix
\section{Proof of Lemma \ref{recur}}\label{A}
By using definition relation \eqref{eq:hermite_numbers},  it is easily calculated that 
\begin{equation}\label{eq:Hk0}
H_{n,0}=H_{n+1,0}=1.
\end{equation}
Since, by convention, $H_{n,-1}=0$, Eq.~\eqref{eq:Hk0} coincides with recurrence relation \eqref{H} for $k=0$.  For $k=1,\ldots,\lfloor n/2\rfloor$, we have
\begin{align}\label{eq:Hk}
(n-2k+2)H_{n,k-1}+H_{n,k}&=\frac{2kn!+n!(n-2k+1)}{2^kk!(n-2k+1)!}=\nonumber\\&=\frac{(n+1)!}{2^kk!(n+1-2k)!}=H_{n+1,k}.
\end{align}
Eqs.~\eqref{eq:Hk0} and \eqref{eq:Hk} constitute the proof of recurrence relation \eqref{H} for even $n$ and $k=0,\ldots,\lfloor (n+1)/2\rfloor$, since $\lfloor (n+1)/2\rfloor=\lfloor n/2\rfloor$ for even $n$. For odd $n=2\ell+1$, Eq.~\eqref{H} has to be also proven for $k=\lfloor (n+1)/2\rfloor=\ell+1$:
\begin{equation}\label{eq:Hodd}
H_{2\ell+1,\ell}=\frac{(2\ell+1)!}{2^\ell\ell!}=\frac{(2\ell+2)!}{2^{\ell+1}(\ell+1)!}=H_{2\ell+2,\ell+1}.
\end{equation}
Since, by convention, $H_{2\ell+1,\ell+1}=0$,  we can easily see that Eq.~\eqref{eq:Hodd} coincides with recurrence relation \eqref{H} for $n=2\ell+1$, $k=\ell+1$.  Thus,  the proof of recurrence relation \eqref{H} for both odd and even $n$,  and for $k=0,\ldots,\lfloor (n+1)/2\rfloor$ is completed.
\begin{rem}
Note that $H_{n+1,0}=1$ is the initial condition supplementing recurrence relation \eqref{H}.  Linear recurrence relation \eqref{H} and initial condition \eqref{eq:Hk0}, under the convention that $H_{n,k}=0$ for $k<0$ or $k>n/2$, constitute a complete definition for signless Hermite coefficients $H_{n,k}$, whose unique solution is formula \eqref{eq:hermite_numbers} (see \cite[Sec. 7.2]{Charalambides2002}).
\end{rem}

\section{Formal derivation of Definition \ref{thm:average}}\label{A:average}
First, we recall that, for a deterministic $C^{\infty}(\mathbb{R}\rightarrow\mathbb{R})$ function $g(x)$, its Taylor series expansion around $x_0$ is alternatively expressed via the translation (shift) pseudodifferential operator (see e.g.  \cite[Sec. 1.1]{Glaeske2006}), first introduced by Lagrange:
\begin{equation}\label{eq:taylor}
g(x)=\sum_{m=0}^{\infty}\frac{\hat{x}^m}{m!}\left.\frac{\id^mg(x)}{\id x^m}\right|_{x=x_0}=\left.\exp\left(\hat{x}\frac{\id}{\id x}\right)g(x)\right|_{x=x_0},
\end{equation}
where $\hat{x}=x-x_0$ is the shift argument. Now, we substitute $X$ as the argument of function $g$,  with $X$ being a random variable with non-zero, in general, mean value $\mu$. By choosing $x_0=\mu$, and taking the average in both sides, Eq.~\eqref{eq:taylor} results into
\begin{equation}\label{eq:taylor_averaged}
\mathbb{E}\left[g(X)\right]=\left.\mathbb{E}\left[\exp\left(\hat{X}\frac{\id}{\id x}\right)\right]g(x)\right|_{x=\mu}=\left.M_{\hat{X}}\left(\frac{\id}{\id x}\right)g(x)\right|_{x=\mu}.
\end{equation}
In Eq.~\eqref{eq:taylor_averaged}, $M_{\hat{X}}(u)$ is identified as the moment-generating function of centered variable $\hat{X}:=X-\mu$ \cite[Sec. 5.5]{Papoulis2002}; $M_{\hat{X}}(u)=\mathbb{E}\left[\exp\left(\hat{X}u\right)\right]$.  Due to its resemblance to the shift operator of Eq.~\eqref{eq:taylor}, we shall call pseudodifferential operator $M_{\hat{X}}(\id/\id x)$ the \textit{averaged shift operator}.  For a Gaussian variable $X$ with variance $\sigma^2$, the moment-generating function of $\hat{X}$ takes the form $M_{\hat{X}}(u)=\exp\left(\sigma^2u^2/2\right)$ \cite[example 5.28]{Papoulis2002}, and Eq.~\eqref{eq:taylor_averaged} is specified into
\begin{equation}\label{eq:pseudo2}
\mathbb{E}\left[g(X)\right]=\left.\exp\left(\frac{\sigma^2}{2}\frac{\id^2}{\id x^2}\right)g(x)\right|_{x=\mu}.
\end{equation}
For the zero-mean value case, $\mu=0$, Eq.~\eqref{eq:pseudo2} results in Eq.~\eqref{eq:pseudo}. 
\section{Proof of Lemma \ref{HC}}\label{B}
By taking into consideration that $C(n,\ell;2)=0$ for $\ell<n/2$ \cite[remark 8.8]{Charalambides2002}, recurrence relation \eqref{eq:recurenceC} for generalized factorial coefficients with parameter 2 is rewritten as
\begin{align}\label{eq:recC}
&C(n+1,\ell;2)=(2\ell-n)C(n,\ell;2)+\nonumber\\&+2C(n,\ell-1;2),  \ \ \text{for} \ \ n\in\mathbb{N}, \ \ \ell=\lceil(n+1)/2\rceil,\ldots,n+1,
\end{align}
supplemented by the final condition deduced from Eq.~\eqref{eq:final_cond}:
\begin{equation}\label{conditionC}
C(n+1,n+1;2)=2^{n+1},  \ \ \text{for} \ \ n\in\mathbb{N}. 
\end{equation}
By performing the index change $k=n+1-\ell$,  Eq.~\eqref{eq:recC} is expressed as
\begin{align}\label{eq:recC2}
C(n+1,n&+1-k;2)=(n-2k+2)C(n,n-(k-1);2)+\nonumber\\&+2C(n,n-k;2),   \ \ k=0,\ldots,\lfloor (n+1)/2\rfloor.
\end{align}
Dividing both sides of Eq.~\eqref{eq:recC2} by $2^{n-k+1}$ results in
\begin{align}\label{eq:recC3}
&\frac{C(n+1,n+1-k;2)}{2^{n+1-k}}=(n-2k+2)\frac{C(n,n-(k-1);2)}{2^{n-(k-1)}}+\nonumber\\&+\frac{C(n,n-k;2)}{2^{n-k}},   \ \ k=0,\ldots,\lfloor (n+1)/2\rfloor,
\end{align}
and dividing Eq.~\eqref{conditionC} by $2^{n+1}$ results in  
\begin{equation}\label{conditionC2}
\frac{C(n+1,n+1;2)}{2^{n+1}}=1,  \ \ \text{for} \ \ n\in\mathbb{N}. 
\end{equation}
Thus, recurrence relation \eqref{eq:recC3} and initial condition \eqref{conditionC2} for $C(n,n-k;2)/2^{n-k}$ coincides with the recurrence relation \eqref{H} and initial condition \eqref{eq:Hk0} for $H_{n,k}$.  This constitutes the proof of Eq.~\eqref{eq:HC}. 


\begin{thebibliography}{12}

\bibitem{Abramowitz1964}
M.~Abramowitz and I.~A. Stegun.
\newblock {\em {Handbook of Mathematical Functions with Formulas, Graphs, and
  Mathematical Tables}}.
\newblock Dover Publications (1983 reprint), 10th edition, 1964.

\bibitem{Athanassoulis2019a}
G.A. Athanassoulis and K.I. Mamis.
\newblock {Extensions of the Novikov-Furutsu theorem, obtained by using
  Volterra functional calculus}.
\newblock {\em Physica Scripta}, 94(11):115217, 2019.
\newblock \url{https://doi.org/10.1088/1402-4896/ab10b5}

\bibitem{Cameron1994}
P.~J. Cameron.
\newblock {\em {Combinatorics: Topics, Techniques, Algorithms}}.
\newblock Cambridge University Press, Cambridge UK, 1994.

\bibitem{Charalambides2002}
Ch.~A. Charalambides.
\newblock {\em {Enumerative combinatorics}}.
\newblock Chapman {\&} Hall/CRC, Boca Raton, 2002.

\bibitem{Glaeske2006}
H.-J. Glaeske, A.P. Prudnikov, and K.A. Sk{\`{o}}rnik.
\newblock {\em {Operational Calculus and Related Topics}}.
\newblock Chapman {\&} Hall/CRC, Boca Raton, 2006.

\bibitem{Jia2017}
D.~Jia, E.~Tang, and A.~Kempf.
\newblock {Integration by differentiation: new proofs, methods and examples}.
\newblock {\em Journal of Physics A: Mathematical and Theoretical}, 50:235201,
  2017.
\newblock \url{https://doi.org/10.1088/1751-8121/aa6f32}  

\bibitem{Klyatskin2005}
V.~I. Klyatskin.
\newblock {\em {Stochastic Equations through the eye of the Physicist}}.
\newblock Elsevier, Amsterdam, 2005.

\bibitem{Mukhopadhyay2021}
N.~Mukhopadhyay.
\newblock {On Rereading Stein's Lemma: Its Intrinsic Connection with
  Cram{\'{e}}r-Rao Identity and Some New Identities}.
\newblock {\em Methodology and Computing in Applied Probability}, 23:355--367,
  2021.
\newblock \url{https://doi.org/10.1007/s11009-020-09830-w}  
  
\bibitem{OEIS}
OEIS Foundation Inc. 
\newblock{Triangle of Bessel numbers read by rows,  Entry A100861 in The On-Line Encyclopedia of Integer Sequences},
\newblock{url: \url{https://oeis.org/A100861}}.  

\bibitem{Papoulis2002}
A.~Papoulis and S.~U. Pillai.
\newblock {\em {Probability, Random variables, and Stochastic processes}}.
\newblock McGraw-Hill, New York, 4th edition, 2002.

\bibitem{Skaltsas2021}
D.~Skaltsas,  G. N.~ Rossopoulos and Ch. I. ~Papadopoulos.
\newblock {A Comparative Study of the Reynolds Equation Solution for Slider and Journal Bearings
with Stochastic Roughness on the Stator and the Rotor}.
\newblock {\em Tribology International}, 167:107410, 2022.
\newblock \url{https://doi.org/10.1016/j.triboint.2021.107410}

\bibitem{Stein1981}
Ch.~M. Stein.
\newblock {Estimation of the Mean of a Multivariate Normal Distribution}.
\newblock {\em The Annals of Statistics}, 9(6):1135--1151, 1981.
\newblock \url{https://doi.org/10.1214/aos/1176345632}

\end{thebibliography}
\end{document}